\documentclass[11pt,a5paper]{article}
\usepackage{ProTex}
\AlProTex{red,<<<>>>,title,list,`}

\usepackage{euler,amsmath,defns,ajr,parskip}
\allowdisplaybreaks

\newcommand{\res}{\text{Res}}
\newcommand{\cL}{{\cal L}}

\newcommand{\Z}[1]{\ensuremath{e^{#1 t}\star}}
\newcommand{\Za}[1]{\ifcase#1 undef0
  \or undef1
  \or\Z{-\rat{27}{10}}
  \or\Z{-\rat{38}{5}}
  \else undefinf \fi}
\newcommand{\Zb}[1]{\ifcase#1 undef0
  \or\Z{-\rat{2}{\epsilon}}
  \or\Z{-\rat{5}{\epsilon}}
  \or\Z{-\rat{10}{\epsilon}}
  \else undefinf \fi}

\title{Computer algebra compares the stochastic superslow manifold of an averaged SPDE with that of the original slow-fast SPDE}

\author{A.~J. Roberts\thanks{School of Mathematical Sciences, University of Adelaide, South Australia.  
\protect\url{http://www.maths.adelaide.edu.au/anthony.roberts}}}

\date{March 22, 2011}

\begin{document}

\maketitle

\begin{abstract}
The computer algebra routines\footnote{The computer algebra routines are written in the free package \textsc{reduce}.  At the time of writing, \url{http://www.reduce-algebra.com/} provides information about \textsc{reduce}.} documented here empower you to reproduce and check many of the details described by an article on large deviations for slow-fast stochastic systems~\cite{Wang2010b}.  We consider a `small' spatial domain with two coupled concentration fields, one governed by a `slow' reaction-diffusion equation and one governed by a stochastic `fast' linear equation.  In the regime of a stochastic bifurcation, we derive two superslow models of the dynamics: the first is of the averaged model of the slow dynamics derived via large deviation principles; and the second is of the original fast-slow dynamics.  Comparing the two superslow models validates the averaging in the large deviation principle in this parameter regime~\cite{Wang2010b}.
\end{abstract}

\tableofcontents

\section{Iterative computer algebra derives the model}
\label{sec:ca}

Construct a one element model of the `slow' stochastic reaction diffusion
equation
\begin{eqnarray}&&
    {\bar u}_t=\partial_{xx} {\bar u} +\lambda\sin {\bar u} -(1-\partial_{xx})^{-1}{\bar u} -\sqrt\epsilon\sigma(1-\partial_{xx})^{-1}\phi( x ,t)
    \label{eq:srd}\\&&
    \mbox{such that}\quad  {\bar u}=0\mbox{ at } x =0,\pi\,,
    \nonumber
\end{eqnarray}
near the deterministic bifurcation that occurs at $\lambda=3/2$\,, to effects quadratic in the noise amplitude~$\sigma$.  We seek the normal form where the evolution involves no convolutions~\cite{Chao95, Roberts06k}.  Also, transform the quadratic noise in the evolution.  Throughout we adopt the Stratonovich interpretation of stochastic differential equations so that the ordinary rules of calculus apply.

The stochastic slow model appears to be, when parameter $\lambda=\rat32+\lambda'$ and upon truncating the noise to just the first three sine modes,
\begin{align}
\dot {\bar a}={}&\lambda'{\bar a}-(\rat3{16}+\rat18\lambda'){\bar a}^3
+\rat{91}{9728}{\bar a}^5
-\sqrt\epsilon\sigma(\rat12\phi_1+\rat3{1216}{\bar a}^2\phi_3)+\cdots
\label{eq:abar}
\end{align}
when the stochastic slow manifold is
\begin{align*}
{\bar u}={}&
{\bar a}\sin x +\rat5{608}{\bar a}^3\sin3x
\\&{}
-\sqrt\epsilon\sigma\left\{ 
\rat15\sin2x\,\Za2\phi_2 +\rat1{10}\sin3x\,\Za3\phi_3
\right.\\&\left.{}
+\lambda'\left[ \rat15\sin2x\,\Za2\Za2\phi_2 +\rat1{10}\sin3x\,\Za3\Za3\phi_3 \right]
\right\}+\cdots
\end{align*}

In outline, the algorithm iteratively determines the stochastic slow manifold model, then finally transforms to a weak model by replacing quadratic noises by their long time equivalents.  Earlier research~\cite{Roberts96a} explained the centre manifold rationale and the computational effectiveness of this simple algorithm, albeit there restricted to deterministic systems.

\<ssmaveq\><<<
% see cassmaveq.pdf for documentation
`<initialisation`>
`<linear noise effects`>
`<quadratic noise effects`>
sig:=small*sigma;
let { small^6=>0 };
it:=1$
repeat begin
`<compute residual`>
`<update ssm`>
    showtime;
end until res=0 or (it:=it+1)>20;
write gssm:=sub(small=1,g);
`<transform quadratic noise`>
end;
>>>

\subsection{Initialisation}

Trivially improve printing.

\<initialisation\><<<
on div; off allfac; on revpri; 
factor sigma,sin,small;
linelength 65$
>>>

Define the parameter~$\lambda$ to be a small perturbation away from critical.  Scale this with ordering parameter~\verb|small| in order to control truncation in the multiple small parameters.

\<initialisation\><<<
lamb:=3/2+small^2*lam;
>>>

Linearise products of trigonometric functions via \verb|trigsimp|.

Define $\alpha_m$ to be the decay rate of linear modes,
here $\alpha_m=m^2-3/2+1/(m^2+1)$\,, so that the spatial modes decay linearly like
$\sin(mx)e^{-\alpha_mt}$.  

\<initialisation\><<<
procedure alfa(m); (m^2-3/2+1/(m^2+1))$ 
>>>

Define the inverse of the linear operator, $\cL^{-1}\sin(mx) =\sin(mx) /\alpha_m$\,, as the linear operator is $\cL=-3/2-\partial_{xx}+(1-\partial_{xx})^{-1}$\,.  Note: we only define and use this for $m\geq 2$\,.

\<initialisation\><<<
operator uinv; linear uinv;
let uinv(sin(~m*x),xt)=>sin(m*x)/alfa(m);
>>>

Define the linear operator, $(1-\partial_{xx})^{-1}\sin(nx)$.

\<initialisation\><<<
operator  iddi; linear  iddi;
let  { iddi(sin(~n*x),x) => sin(n*x)/(1+n^2)
     , iddi(sin(x),x) => sin(x)/2 };
>>>

Paramterise solutions by an evolving amplitude~$\bar a(t)$ (or `order
parameter').  Its evolution is $d{\bar a}/dt =\dot {\bar a} =g$\,.

\<initialisation\><<<
depend a,t;
let df(a,t)=>g;
>>>

Then the most basic linear approximation to the dynamics on the element
is ${\bar u}={\bar a}\sin x$ where $\dot {\bar a}=0$\,.  Scale the amplitude to be small.

\<initialisation\><<< 
u:=small*a*sin(x); 
v:=u/2; 
g:=0; 
>>>

\subsection{Compute residual}

The parameter~\verb|small| conveniently controls the trunction in nonlinearity and other small parameters.  The iteration terminates when the residual of the reaction diffusion equation is zero to the specified order of smallness.   Note: this parametrisation with~\verb|small| should create deterministic models to errors~$\Ord{{\bar a}^4+\epsilon^2}$, or equivalent, as we scale $\epsilon=\verb|eps|$ with \verb|small^2|. For some strange reason we need to do some operation on \verb|res| in order for relevant terms to cancel---here I use \verb|trigsimp|, but something else might serve.  

\<compute residual\><<<
    sinu:=trigsimp(u-u^3/6+u^5/120-u^7/5040,combine);
    res:=-df(u,t)+df(u,x,2)+lamb*sinu-iddi(u,x)
        -small*rooteps*sig*iddi(noise,x);
    res:=trigsimp(res,combine);
    write lengthres:=length(res);
>>>

Define~$\sqrt\epsilon$ which scales the size of the noise.  Looks like we do not have to worry about \verb|rooteps/eps| not simplifying.  

\<initialisation\><<<
let rooteps^2=>eps;
>>>

\subsection{Update the stochastic slow manifold}
   
Let $T$, \verb|tt|, label the fast time of stochastic fluctuations so we can separate the stochastic fluctuations from the superslow evolution of the amplitude~$\bar a$.  Also introduce \verb|xt| to label both the subgrid spatial scales and  time scales so we can group all factors in the space-time dynamics.

\<initialisation\><<< 
depend tt,t;
depend x,xt; 
>>>
   
Then update driven by the residual.  Divide the evolution~\verb|g| by \verb|small| to best keep track of the correct counting of the `order' of a term.

\<update ssm\><<<
    g:=g+(gd:=secular(res,xt))/small;
    u:=u+uinv(res-gd*sin(x),xt);
>>>

\subsection{Linear noise effects}

Introduce the noise in its spatial Fourier decomposition 
\begin{displaymath}
    \phi(x,t)=\sum_{n=1}^\infty \phi_n(t)\sin nx\,.
\end{displaymath}
Parametrise the amplitude with~$\sigma$.   Truncate the spatial structure of noise. Three terms in the noise appears adequate to show the typical interactions between noise and other dynamics.  Have not explored details of better resolution of the noise.  

\<linear noise effects\><<<
operator phi; depend phi,tt,xt;
noise:=for n:=1:3 sum phi(n,{})*sin(n*x);
>>>

Let  \verb|phi(n,{m1,...})|
denote convolutions with \verb|exp(-m1*t)...|, that is,
\begin{displaymath}
    \phi_{n,(m_1,m_2,\ldots)}
    =\exp(-{m_1}t)\star \exp(-{m_2}t)\star 
    \cdots\star\phi_n(t)\,;
\end{displaymath}
so 
\begin{displaymath}
    \partial_t \phi_{n,(m_1,m_2,\ldots)}
    =-{m_1} \phi_{n,(m_1,m_2,\ldots)} +\phi_{n,(m_2,\ldots)}\,.
\end{displaymath}

But if we pull out a decay rate in~$1/\epsilon$ then keep bookkeeping correct by dividing by~\verb|small^2| unless it is a \emph{sole} convolution by~\Ord{1/\epsilon}.  This latter case is only used in the next section.

\<linear noise effects\><<<
let { df(phi(~m,~p),t)=>df(phi(~m,~p),tt)
    , df(phi(~m,~p),tt)=>(-first(p)*phi(m,p)+phi(m,rest(p)))
        when deg(1/first(p),eps)=0 
    , df(phi(~m,~p),tt)=>(-first(p)*phi(m,p)/small^2
        +phi(m,rest(p))/small^(if rest(p)={} then 1 else 2))
        when deg(1/first(p),eps)=1
    };
>>>

Recall the equation for updates ${\bar u}'$~and~$g'$ is $g' +\cL{{\bar u}'} ={}$residual, where now the operator $\cL=\partial_t-3/2-\partial_{xx}+(1-\partial_{xx})^{-1}$ includes fast time variations.  The operator \verb|secular| extracts from the residual all those terms which would generate generate secular growth in the field~$u$ and so instead must be placed in the model's evolution~$g$.  The last rule here comes from integration by parts and is essential in order to eliminate memory integrals (convolutions) in the model evolution.

The if-clause in the last is needed for the next section to account for the convolution on the fast time only being of~\Ord{\sqrt\epsilon}, and so order increases

\<linear noise effects\><<<
operator secular; linear secular;
let { secular(sin(~m*x),xt)=>0
    , secular(sin(~m*x)*~aa,xt)=>0
    , secular(sin(x),xt)=>1
    , secular(sin(x)*phi(~n,~p),xt)=>
      phi(n,{})*(for each r in p product (1/r))
      *(if p neq{} and deg(1/first(p),eps)=1 then small else 1)
    };
>>>

Extend the inverse operator to terms with fast time variations as well as fast (subgrid) space variations.  Recursive procedure \verb|gungb| extracts the non-secular parts of fluctuating $\sin x$ components. Have to adjust the smallness whenever the convolution transformed is on the $\epsilon$~scale.

\<linear noise effects\><<<
procedure gungb(n,p);
    if p={} then 0 else 
    (gungb(n,rest(p))-phi(n,p)
    *(if deg(1/first(p),eps)=0 then 1 else small^2)
    )/first(p)$
    
let { uinv(sin(~m*x)*phi(~n,~p),xt)=>phi(n,(alfa(m)).p)*sin(m*x)
    , uinv(sin(x)*phi(~n,~p),xt)=>gungb(n,p)*sin(x)
    };
>>>

\subsection{Quadratic noise effects}

Now let $Z_p$ denote multiple convolutions of in time of any term,
\verb|zz(a,p)| (though I only use~$Z$ for quadratic terms, it may well
be able to replace the linear convolutions).  That is,
\begin{displaymath}
    Z_{(m_1,m_2,\ldots)} =\exp(-{m_1}t)\star Z_{(m_2,\ldots)}
    \quad\text{and}\quad
    Z_{(\,)}=1\,.
\end{displaymath}

\<quadratic noise effects\><<<
operator zz; depend zz,tt,xt; 
let { zz(~a,{})=>a
    , df(zz(~a,~p),t)=>df(zz(a,p),tt)
    , df(zz(~a,~p),tt)=>-first(p)*zz(a,p)+zz(a,rest(p))
        when deg(1/first(p),eps)=0 
    };
>>>

To extract quadratic corrections to the evolution, use integration by
parts so all non-integrable convolutions are reduced to the cannonical
form of the convolution being entirely over one noise in a quadratic
term, either $\phi_n\phi_{m,(\ldots)}$ or $\phi_{n,(\ldots)}\phi_m$.

Have now made this very complicated for at least some of the cases when the convolutions may be over $\epsilon$-fast time scales. 

\<quadratic noise effects\><<<
procedure gungd(n,p,m,q); 
    if (p={})or(q={}) then phi(n,p)*phi(m,q) 
    else if deg(1/first(p),eps)=deg(1/first(q),eps) then
    (gungd(n,rest(p),m,q)+gungd(n,p,m,rest(q)))
    /(first(p)+first(q))
    *(if deg(1/first(p),eps)=0 then 1 else small^2)
    else if deg(1/first(p),eps)=1 then
    (gungd(n,rest(p),m,q)*(if rest(p)={} then small else small^2)
      +gungd(n,p,m,rest(q))*small^2
      )/first(p)*sub(rat=-first(q)/first(p),geom)
    else gungd(m,q,n,p)$
    
let { secular(sin(x)*zz(~a,~p),xt) =>secular(sin(x)*a,xt)
      *(for each r in p product (1/r))
    , secular(sin(~m*x)*zz(~a,~p),xt)=>0
    , secular(sin(x)*phi(~n,~p)*phi(~m,~q),xt) =>gungd(n,p,m,q) 
    , secular(sin(x)*phi(~n,~p)^2,xt) =>gungd(n,p,n,p)
    };
>>>

Extend $\cL^{-1}$ operator \verb|uinv| to handle quadratic terms.
First, integration by parts gives all integrable contributions from
direct product terms.

Have to similarly modify \verb|gunge| for at least some of the cases when the convolutions may be over $\epsilon$-fast time scales. 

\<quadratic noise effects\><<<
procedure gunge(n,p,m,q); 
    if (p={})or(q={}) then 0 
    else  if deg(1/first(p),eps)=deg(1/first(q),eps) then
    (-phi(n,p)*phi(m,q)
      *(if deg(1/first(p),eps)=0 then 1 else small)^2
      +gunge(n,rest(p),m,q)
      +gunge(n,p,m,rest(q))
    )/(first(p)+first(q))
    else if deg(1/first(p),eps)=1 then
    (-phi(n,p)*phi(m,q)*small^2
      +gunge(n,rest(p),m,q)*(if rest(p)={} then small else small^2)
      +gunge(n,p,m,rest(q))*small^2
      )/first(p)*sub(rat=-first(q)/first(p),geom)
    else gunge(m,q,n,p)$
let { uinv(sin(x)*phi(~n,~p)*phi(~m,~q),xt)
      =>gunge(n,p,m,q)*sin(x)
    , uinv(sin(x)*phi(~n,~p)^2,xt)=>gunge(n,p,n,p)*sin(x)
    };
>>>

Second, similar integration by parts gives integrable contribution
from terms involving convolutions of products.

\<quadratic noise effects\><<<
procedure gungf(a,p);
    if p={} then 0 else 
    (gungf(a,rest(p))-zz(a,p))/first(p)$
let { uinv(sin(~l*x)*phi(~n,~p)*phi(~m,~q),xt)
      =>sin(l*x)*zz(phi(n,p)*phi(m,q),{alfa(l)})
    , uinv(sin(~l*x)*phi(~n,~p)^2,xt)
      =>sin(l*x)*zz(phi(n,p)^2,{alfa(l)})
    , uinv(sin(~l*x)*zz(~a,~p),xt)=>sin(l*x)*zz(a,alfa(l).p)
    , uinv(sin(x)*zz(~a,~p),xt)=>sin(x)*gungf(a,p)
    };
>>>

\subsection{Transform quadratic noise}
\label{sec:transquad}

Now proceed to transform the strong model to a weak model by replacing
the quadratic noises by their effective long term drift and volatility.
Earlier work~\cite{Roberts05e} determines the rationale for the details of this
transformation.

Set $\verb|small|=1$ as it has done its job of truncating the nonlinear
terms in the asymptotic expansion.

\<transform quadratic noise\><<<
small:=1;
write "Now transforming the quadratic noises";
>>>

Now transform the quadratic noise into new noises~$\psi$ (\verb|psi|)
that are equivalent in their long time statistics: the
operator~\verb|long| implements the long-time equivalent noises as
determined earlier~\cite{Roberts05e}.  For now only transform up to two
convolutions.  These new noises have subscripts that uniquely identify
them.

\<transform quadratic noise\><<<
operator long; linear long;
operator psi; depend psi,tt,xt;
let { long(1,tt)=>1
    , long(phi(~i,{}),tt)=>phi(i,{})
    , long(phi(~i,{})*phi(~j,{~k}),tt)
      => 1/2*(if i=j then 1 else 0)
      +psi(i,j,{k})/sqrt(2*k)
    , long(phi(~i,{})*phi(~j,{~k2,~k1}),tt)
      => (psi(i,j,{k1})/sqrt(2*k1)
      +psi(i,j,{k2,k1})/sqrt(2*k2))/(k1+k2)
    };
gg:=long(g,tt)$
>>>

Root sum squares of the determined noise coefficients; this procedure implicitly assumes that there is no correlation between the multitude of noises in these two terms in the amplitude equation.  This assumption appears correct for these two terms in this \pde.  In general we would need to do a $QR$~factorisation of the noise terms.

\<transform quadratic noise\><<<
operator sumsqpsi; linear sumsqpsi;
let { sumsqpsi(1,tt)=>0
    , sumsqpsi(psi(~i,~j,~p),tt)=>0
    , sumsqpsi(psi(~i,~j,~p)^2,tt)=>1
    , sumsqpsi(psi(~i,~j,~p)*psi(~ii,~jj,~pp),tt)=>0
    };
>>>

Have a look at the numerical coefficients.

\<transform quadratic noise\><<<
on rounded; print_precision 5;
gg:=gg;
>>>

Extract the coefficients of the terms in $\sigma^2$~and~$\sigma^2a$, 
both mean and fluctuating.

\<transform quadratic noise\><<<
let abs(eps)=>eps;
c20:=sqrt(sumsqpsi(coeffn(coeffn(gg,sig,2),a,0)^2,tt));
c21mean:=(coeffn(coeffn(gg,sig,2),a,1) 
    where psi(~i,~j,~p)=>0);
c21:=sqrt(sumsqpsi(coeffn(coeffn(gg,sig,2),a,1)^2,tt));
>>>

Switch back to the rational arithmetic mode for any other postprocessing.

\<transform quadratic noise\><<<
off rounded;
showtime;
>>>

%%%%%%%%%%%%%%%%%%%%%%
\OutputCode\<ssmaveq\>
%%%%%%%%%%%%%%%%%%%%%%

Executing the resultant code constructs the superslow model of the stochastic bifurcation in the slow averaged \spde{}s.

\section{Model interacting fast-slow-superslow components}

This section constructs a superslow model of the fast-slow stochastic reaction diffusion
equation
\begin{align}&
    u_t=\partial_{xx} u +\lambda\sin u -v\,,
    \label{eq:srdu}\\&
    \epsilon v_t=\partial_{xx} v -v +u +\sqrt\epsilon\sigma\phi( x ,t),
    \label{eq:srdv}\\&
    \mbox{such that}\quad  u=v=0\mbox{ at } x =0,\pi\,,
    \nonumber
\end{align}
near the deterministic bifurcation that occurs at $\lambda=3/2$\,, to effects quadratic in the noise amplitude~$\sigma$, and seeks the normal form where the evolution involves no convolutions.  Recall that throughout we adopt the Stratonovich interpretation of stochastic differential equations so that the ordinary rules of calculus apply.

The presence of the three time scales in the dynamics---the fast~$v$, the slow~$u$, and the superslow evolution of the bifurcation amplitude~$a$---means that the algorithm outlined here is one of the most technically challenging ones I have implemented.  I recommend understanding simpler systems before attempting to understand the details of this section.

The resulting model appears to be the following to some order in small parameters.  In terms of the superslow evolving amplitude~$a(t)$, where $u\approx a\sin x$ and $v\approx \rat12 a\sin x$\,, parameter $\lambda=\rat32+\lambda'$, and noise in just three sine modes,  a stochastic differential equation for the amplitude is
\begin{align}
\dot a={}&\lambda'(1+\rat14\epsilon\lambda')a 
-(\rat3{16}+\rat18\lambda'+\rat3{64}\epsilon)a^3
+\rat{91}{9728}a^5
\nonumber\\&{}
-\sqrt\epsilon\sigma(\rat12+\rat18\epsilon
-\rat14\epsilon\lambda' +\rat9{64}\epsilon a^2 %+\rat1{32}\epsilon^2
)\phi_1
\nonumber\\&{}
-\sqrt\epsilon\sigma(\rat3{1216}+\rat3{4864}\epsilon)a^2\phi_3 
\nonumber\\&{}
+\epsilon\sigma^2a\left[-\rat1{180}\phi_2\Za2\phi_2
+\rat3{1216}\phi_1\Za3\phi_3
-\rat6{6080}\phi_3\Za3\phi_3
\right]
\nonumber\\&{}
+\cdots
\end{align}
To errors~\Ord{\epsilon}, this evolution equation is identical to the slow model~\eqref{eq:abar} of the averaged equation with fluctuations.

The corresponding stochastic superslow manifold appears to be that the slow field
\begin{align*}
 u={}& a\sin x 
+\rat5{608}a^3\sin 3x
+\rat12\sqrt\epsilon\sigma\sin x\,\Zb1\phi_1
\\&{}
-\rat15\sqrt\epsilon\sigma\sin2x\left[\Za2-\Zb2\right]\phi_2
\\&{}
-\rat1{10}\sqrt\epsilon\sigma\sin3x\left[\Za3-\Zb3\right]\phi_3 
+\cdots
\end{align*}
whereas the fast field has \Ord1~fluctuations 
\begin{align*}
v={}&\rat12a\sin x
+\rat1{1216}a^3\sin3x
\\&{}
+\frac\sigma{\sqrt\epsilon}\sin x\left[(1+\rat14\epsilon)\Zb1 +\rat12\Zb1\Zb1\right]\phi_1
\\&{}
+\frac\sigma{\sqrt\epsilon}
\sin2x\left[(1+\rat1{25}\epsilon)\Zb2 -\epsilon\rat1{25}\Za2 +\rat15\Zb2\Zb2 \right]\phi_2
\\&{}
+\frac\sigma{\sqrt\epsilon}
\sin3x\left[(1+\rat1{100}\epsilon)\Zb3 -\epsilon\rat1{100}\Za3 +\rat1{10}\Zb3\Zb3 \right]\phi_3
\\&{}
+\cdots
\end{align*}

In outline, the algorithm iteratively determines the stochastic superslow manifold model~\cite[e.g.]{Roberts96a}.

Seem to have to keep one or two orders higher in \verb|small| as there is a division somewhere.  So actually compute to residuals and errors one or two orders in \verb|small| less than apparently allowed here---this fudges the computations so that enough is kept to do the cancellation, then later truncates. That is, \verb|small^6=>0| actually computes to errors~$\Ord{\texttt{small}^4}$.  Here, for some reason we have to tread carefully to get up to fifth order terms in~\verb|small|: first, compute linear noise effects; then, second, when residuals are zero, up the order to retain quadratic terms in noise and continue iterating.  Alternatively, we could just seek terms up to fourth order in~\verb|small| by \verb|let small^7=>0|.  But I do want to get to fifth order because of the challenge. 

\<ssmuv\><<<
% see cassmaveq.pdf for documentation
let { sigma^2=>0, small^8=>0 };
sig:=small*sigma;
`<initialisation`>
`<linear noise effects`>
`<quadratic noise effects`>
it:=1$
repeat begin
`<update from fast residual`>
if {resu,resv}={0,0} then clear sigma^2;%implicitly sigma^3=>0;
`<update from slow residual`>
    showtime;
end until {resu,resv}={0,0} and (sigma^2neq0) or (it:=it+1)>19;
%write ussm:=sub(small=1,u);
%write vssm:=sub(small=1,v);
write gssm:=sub(small=1,g);
`<transform quadratic noise`>
end;
>>>

\subsection{Some initialisation things}

Define $\beta_m$ to be the relative decay rate of linear modes of the fast variable~$v$ on the element, here $\beta_m=m^2+1$\,, so that the spatial modes in~$v$ decay linearly like $\sin(mx)\exp(-\beta_mt/\epsilon)$\,.  

\<initialisation\><<<
procedure beta(m); (m^2+1)$
>>>

Define some of the inverse of a linear operator.  Now there are significant subtleties here: each convolution with rates~\Ord{1/\epsilon} are themselves. Thus smallness is hidden in the convolution rates; consequently we have to track them artificially though a parameter such as~\verb|small|.  Use~\verb|small| to count both the direct~$\epsilon$ and the hidden ones in the convolutions, as well as the other small parameters.  

\begin{table}
\caption{order of magnitude of convolution operators.}
\label{tbl:omco}
\begin{center}
\begin{tabular}{l|cccc}
& 1 & \Z{-\alpha} & $(\Z{-\alpha})^2$ & $(\Z{-\alpha})^3$ \\
\hline
1 & 1 & 1 & 1 & 1 \\
\Z{-\frac\beta\epsilon} & \Ord{\epsilon^{1/2}} & \Ord{\epsilon} & \Ord{\epsilon} & \Ord{\epsilon} \\
$(\Z{-\frac\beta\epsilon})^2$ & \Ord{\epsilon^{3/2}} & \Ord{\epsilon^2} & \Ord{\epsilon^2} & \Ord{\epsilon^2} \\
$(\Z{-\frac\beta\epsilon})^3$ & \Ord{\epsilon^{5/2}} & \Ord{\epsilon^3} & \Ord{\epsilon^3} & \Ord{\epsilon^3} \\
\hline
\end{tabular}
\end{center}
\end{table}

A further complication is that single bare convolution is actually~\Ord{\sqrt\epsilon} \cite[equation~(27)]{Roberts05c}; Table~\ref{tbl:omco} lists the correct order of magnitude of various convolutions.  This complication is simplified a little by separating convolutions that occur on the fast time from the convolutions that occur on the low time scale.

The linear equations for updates are
\begin{align*}
-g\sin x-u_t&{}+u_{xx}+\rat32 u-v+\res_u=0\,, \\
-\epsilon v_t&{}+u+v_{xx}-v+\res_v=0\,.
\end{align*}
When considering mode~$\sin mx$, the linear equations for updates are, including a correction~$g$ to the evolution only in the case of the critical $m=1$\,,
\begin{align*}
-g-u_t&{}-(\alpha_m-1/\beta_m)u-v+\res_u=0\,, \\
-\epsilon v_t&{}+u-\beta_m v+\res_v=0\,,
\end{align*}
for the previous defined constants $\alpha_m$~and~$\beta_m$. The details of solving for updates are not critical to correctness of the results (as the results should only depend upon driving the residuals to zero), but the details will determine whether the iteration does converge to zero the residuals.

\paragraph{Updates from the $u$~equation}
For residuals of the slow-equation, make updates to the $u$-field driven by the $u$-residual, and correspondingly update the $v$-field in a way that will not change its $v$-residual at this order.  We insist on not changing the $v$-residual because this update is considered second and we must not undo earlier corrections driven from the $v$-residual.  Dividing the $u$-update by~$\beta_m$ is sufficient for the $v$-update.  List here first the deterministic updates, second the generic linear noise update, and last the updates for resonant terms.  Use procedure \verb|gungb| to  extract the non-resonant parts of~$\res_u$.

Account for smoothing effect of convolution on noise via the if-clauses.  Table~\ref{tbl:omco} shows that when a term goes from multiple fast time convolutions to include one additional slow time convolution, then the order of the term increases by~$\sqrt\epsilon$ (the following provision assumes we linearise convolutions so that any one term only has convolutions on the same time scale). 

\<linear noise effects\><<<
operator uuinv; linear uuinv;
operator vuinv; linear vuinv;
let { uuinv(sin(~m*x),xt)=>sin(m*x)/alfa(m)
    , vuinv(sin(~m*x),xt)=>sin(m*x)/alfa(m)/beta(m)
    , uuinv(sin(~m*x)*phi(~n,~p),xt)
        => phi(n,(alfa(m)).p)*sin(m*x)
        *(if p neq{} and deg(1/first(p),eps)=1 then small else 1)
    , vuinv(sin(~m*x)*phi(~n,~p),xt)
        => phi(n,(alfa(m)).p)*sin(m*x)/beta(m)
        *(if p neq{} and deg(1/first(p),eps)=1 then small else 1)
    , uuinv(sin(x)*phi(~n,~p),xt)=>gungb(n,p)*sin(x)
    , vuinv(sin(x)*phi(~n,~p),xt)=>gungb(n,p)*sin(x)/beta(1)
    };
>>>

To deal with quadratic noise, the following appear to to be enough for errors no higher order than~\verb|small|${}^8$.   That is, with the modifications made to \verb|gunge| and \verb|gungd|.  We do not seem to need any extra transformations from the residual of the fast $v$-equation, probably because it is linear.

\<quadratic noise effects\><<<
let { uuinv(sin(~m*x)*phi(~n,~p)*phi(~l,~q),xt)
        => zz(phi(n,p)*phi(l,q),{alfa(m)})*sin(m*x)
        %*(if p neq{} and deg(1/first(p),eps)=1 then small else 1)
    , vuinv(sin(~m*x)*phi(~n,~p)*phi(~l,~q),xt)
        => zz(phi(n,p)*phi(l,q),{alfa(m)})*sin(m*x)/beta(m)
        %*(if p neq{} and deg(1/first(p),eps)=1 then small else 1)        
    , uuinv(sin(~m*x)*zz(~n,~p),xt) => zz(n,alfa(m).p)*sin(m*x)
    , vuinv(sin(~m*x)*zz(~n,~p),xt) 
      => zz(n,alfa(m).p)*sin(m*x)/beta(m)
    , uuinv(sin(~m*x)*phi(~n,~p)^2,xt)
        => zz(phi(n,p)^2,{alfa(m)})*sin(m*x)
        %*(if p neq{} and deg(1/first(p),eps)=1 then small else 1)
    , vuinv(sin(~m*x)*phi(~n,~p)^2,xt)
        => zz(phi(n,p)^2,{alfa(m)})*sin(m*x)/beta(m)
        %*(if p neq{} and deg(1/first(p),eps)=1 then small else 1)
    , uuinv(sin(x)*phi(~n,~p)*phi(~l,~q),xt)=>gunge(n,p,l,q)*sin(x)
    , vuinv(sin(x)*phi(~n,~p)*phi(~l,~q),xt)
      => gunge(n,p,l,q)*sin(x)/beta(1)
    , uuinv(sin(x)*phi(~n,~p)^2,xt)=>gunge(n,p,n,p)*sin(x)
    , vuinv(sin(x)*phi(~n,~p)^2,xt)=>gunge(n,p,n,p)*sin(x)/beta(1)
    , uuinv(sin(x)*zz(~n,~p),xt)=> 
        ( uuinv(sin(x)*zz(n,rest(p)),xt)
          -zz(n,p)*sin(x) )/first(p)
    , vuinv(sin(x)*zz(~n,~p),xt)=> 
        ( vuinv(sin(x)*zz(n,rest(p)),xt)
          -zz(n,p)*sin(x) )/first(p)/beta(1)
    };
>>>

\paragraph{Updates from the $v$~equation}
Corrections to the $u$~and~$v$ fields arise from the $v$-residual.  However, because we consider this residual first in each iteration (not that first makes a lot of sense in  an iterative loop), we are free to modify field~$u$ in a way that would affect the residuals at the same order.  The key aspect is that we must not affect the residuals at a lower order in the $u$-residual---achieving this aspect is hard enough, which is why I implement corrections from the $v$-residual first.

First define the deterministic updates.

\<linear noise effects\><<<
operator vvinv; linear vvinv;
operator uvinv; linear uvinv;
let { vvinv(sin(~m*x),xt)=>sin(m*x)*(alfa(m)-1/beta(m))
        /alfa(m)/beta(m)
    , uvinv(sin(~m*x),xt)=>-sin(m*x)/alfa(m)/beta(m)
    , vvinv(sin(x),xt)=>sin(x)/beta(1)
    , uvinv(sin(x),xt)=>0
>>>

Second deal with the variety of linear noise terms.  When a noise term in the residual is a convolution over the slow-scale, then the convolution is smooth and its time derivative correspondingly of the same order so that the $\epsilon v_t$ causes no problem.

\<linear noise effects\><<<
    , vvinv(sin(~m*x)*phi(~n,~p),xt) 
        => phi(n,p)*sin(m*x)/beta(m) 
        when p neq {} and deg(1/first(p),eps)=0
    , uvinv(sin(~m*x)*phi(~n,~p),xt) => 0 
        when p neq {} and deg(1/first(p),eps)=0
>>>

The critical mode is no different when forcing in the $v$-residual.

\<linear noise effects\><<<
    , vvinv(sin(x)*phi(~n,~p),xt)
        => phi(n,p)*sin(x)/beta(1) 
        when p neq {} and deg(1/first(p),eps)=0
    , uvinv(sin(x)*phi(~n,~p),xt) => 0 
        when p neq {} and deg(1/first(p),eps)=0
>>>

However, when the noise term in the $v$-residual is not smooth, either because it is a bare white noise or because it is a convolution over the fast time scale, then we must be more careful because the $\epsilon v_t$~term is important.  Because the update has to be relatively large,\footnote{I conjecture that it is this that affects the management of the smallness parameter.} we have to use the $u$-field to cancel the effect in the $u$-residual of updates from the $v$~equation.

Here if the convolution is the first fast-time convolution, \verb|p={}|, then choose the correct scale in \verb|small| as then the convolution is only~\Ord{\sqrt\epsilon}.

\<linear noise effects\><<<
    , vvinv(sin(~m*x)*phi(~n,~p),xt)
        => phi(n,(beta(m)/eps).p)*sin(m*x)/eps
        /(if p={} then small else 1)
        when p={} or deg(1/first(p),eps)=1
    , uvinv(sin(~m*x)*phi(~n,~p),xt)
        => phi(n,(beta(m)/eps).p)*sin(m*x)/beta(m)
        *(if p={} then small else small^2)
        when p={} or deg(1/first(p),eps)=1
>>>

The critical mode is no different when forcing in the $v$-residual.

\<linear noise effects\><<<
    , vvinv(sin(x)*phi(~n,~p),xt)
        => phi(n,(beta(1)/eps).p)*sin(x)/eps
        /(if p={} then small else 1)
        when p={} or deg(1/first(p),eps)=1
    , uvinv(sin(x)*phi(~n,~p),xt)
        => phi(n,(beta(1)/eps).p)*sin(x)/beta(1)
        *(if p={} then small else small^2)
        when p={} or deg(1/first(p),eps)=1
    };
>>>

Second deal with a variety of quadratic noise terms.  When a noise term in the residual is a convolution over the slow-scale, then the convolution is smooth and its time derivative correspondingly of the same order so that the $\epsilon v_t$ causes no problem.

\<quadratic noise effects\><<<
let { vvinv(sin(~m*x)*zz(~n,~p),xt) 
        => zz(n,p)*sin(m*x)/beta(m) when deg(1/first(p),eps)=0
    , uvinv(sin(~m*x)*zz(~n,~p),xt) => 0 when deg(1/first(p),eps)=0
    , vvinv(sin(~m*x)*phi(~n,~p)^2,xt) 
        => phi(n,p)^2*sin(m*x)/beta(m) 
        when p neq{} and deg(1/first(p),eps)=0
    , uvinv(sin(~m*x)*phi(~n,~p)^2,xt) => 0 
        when p neq{} and deg(1/first(p),eps)=0
    , vvinv(sin(~m*x)*phi(~n,~p)*phi(~l,~q),xt) 
        => phi(n,p)*phi(l,q)*sin(m*x)/beta(m) 
        when p neq{} and q neq{} 
        and deg(1/first(p),eps)+deg(1/first(q),eps)=0
    , uvinv(sin(~m*x)*phi(~n,~p)*phi(~l,~q),xt) => 0 
        when p neq{} and q neq{} 
        and deg(1/first(p),eps)+deg(1/first(q),eps)=0
>>>

The critical mode is no different when forcing in the $v$-residual.

\<quadratic noise effects\><<<
    , vvinv(sin(x)*zz(~n,~p),xt)
        => zz(n,p)*sin(x)/beta(1) when deg(1/first(p),eps)=0
    , uvinv(sin(x)*zz(~n,~p),xt) => 0 when deg(1/first(p),eps)=0
    , vvinv(sin(x)*phi(~n,~p)^2,xt) 
        => phi(n,p)^2*sin(x)/beta(1) 
        when p neq{} and deg(1/first(p),eps)=0
    , uvinv(sin(x)*phi(~n,~p)^2,xt) => 0 
        when p neq{} and deg(1/first(p),eps)=0
    , vvinv(sin(x)*phi(~n,~p)*phi(~l,~q),xt) 
        => phi(n,p)*phi(l,q)*sin(x)/beta(1) 
        when p neq{} and q neq{} 
        and deg(1/first(p),eps)+deg(1/first(q),eps)=0
    , uvinv(sin(x)*phi(~n,~p)*phi(~l,~q),xt) => 0 
        when p neq{} and q neq{} 
        and deg(1/first(p),eps)+deg(1/first(q),eps)=0
    };
>>>

\paragraph{Critical: linearise convolutions over different time scales}
I contend that we also want to simplify the convolutions because convolutions of the fast time scale~$\epsilon$ are qualitatively different from convolutions over slow time scales.  Thus we do the following `linearisation' of convolutions: whenever the first two convolutions are over different time scales, we transform the convolution into the sum of two convolutions.  Change of variables in integration shows that
\begin{equation*}
\Z{-\alpha}\Z{-\beta}=\frac1{\beta-\alpha}\left[\Z{-\alpha}-\Z{-\beta}\right].
\end{equation*}
I have not used this transform in other applications because of the necessity to avoid division by zero when $\alpha=\beta$\,.\footnote{However,  in general, maybe I should do this linearisation in order to reduce expressions to a more canonical form.}  Here, we are concerned with convolutions over different time scales and so use this formula where, for example, rate~$\beta$ is replaced by fast rate~$\beta/\epsilon$:
\begin{align*}
\Z{-\alpha}\Z{-\frac\beta\epsilon}&{}=
\frac1{\frac\beta\epsilon-\alpha}\left[\Z{-\alpha}-\Z{-\frac\beta\epsilon}\right]
\\&{}=\frac{\epsilon/\beta}{1-\epsilon\alpha/\beta}\left[\Z{-\alpha}-\Z{-\frac\beta\epsilon}\right].
\end{align*}
Thus we need to divide by $1-r$ for various $r=\epsilon\alpha/\beta$ so store its power series in the variable~\verb|geom|, with \verb|small| to account for the powers of~$\epsilon$.  \emph{The correctness of the following transformation is critical.}

\<linear noise effects\><<<
geom:=for n:=0:deg((1+small^2)^9,small)/2 sum (rat*small^2)^n$
let { phi(~n,~p) => (phi(n,first(p).rest(rest(p)))
        -phi(n,second(p).rest(rest(p)))*(if rest(rest(p))={} 
            or deg(1/first(rest(rest(p))),eps)=1 
            then small else small^2)
        )*sub(rat=first(p)/second(p),geom)/second(p)
      when length(p)>1 and deg(1/first(p) ,eps)=0 
                       and deg(1/second(p),eps)=1
    , phi(~n,~p) => (phi(n,second(p).rest(rest(p)))
        -phi(n,first(p).rest(rest(p)))*(if rest(rest(p))={} 
            or deg(1/first(rest(rest(p))),eps)=1 
            then small else small^2)
        )*sub(rat=second(p)/first(p),geom)/first(p)
      when length(p)>1 and deg(1/first(p) ,eps)=1 
                       and deg(1/second(p),eps)=0
    };
>>>

\subsection{Update from residuals of the fast equation}

The parameter~\verb|small|, controls the truncation in nonlinearity and in small parameters.  The iteration terminates when the residual of the reaction diffusion equation is zero to the specified order of nonlinearity.   For some reason we need to do something nontrivial to the residual in order to force cancellation of terms so I apply \verb|trigsimp|.  Note the multiplication and division by \verb|small| in order to cater for other divisions by~\verb|small| affecting the error truncation.

\<update from fast residual\><<<
    resv:=-small^2*eps*df(v,t)+df(v,x,2)-v+u
        +small*rooteps*sig*noise; 
    resv:=trigsimp(small^2*resv,combine)/small^2;
    write lengthresv:=length(resv);
    u:=u+uvinv(resv,xt);
    v:=v+vvinv(resv,xt);
>>>

\subsection{Update from residuals of the slow equation}
   
Similarly update from the residual of the slow equation.  Divide the evolution~\verb|g| by \verb|small| to best keep track of the correct counting of the `order' of a term.

\<update from slow residual\><<<
    sinu:=trigsimp(u-u^3/6+u^5/120-u^7/5040,combine);
    resu:=-df(u,t)+df(u,x,2)+lamb*sinu-v; 
    resu:=trigsimp(small^2*resu,combine)/small^2;
    write lengthresu:=length(resu);
    g:=g+(gd:=secular(small^2*resu,xt)/small^2)/small;
    u:=u+uuinv(resu-gd*sin(x),xt);
    v:=v+vuinv(resu-gd*sin(x),xt);
>>>

Executing the resultant code constructs the superslow model of the stochastic bifurcation in the fast-slow system of \spde{}s.

%%%%%%%%%%%%%%%%%%%%
\OutputCode\<ssmuv\>
%%%%%%%%%%%%%%%%%%%%

\paragraph{Acknowledgement} This research is supported by Australian Research Council grants DP0774311 and DP0988738.

%%%%%%%%%%%%%%%%%%%%%%%%%%%%

\bibliographystyle{plain}
\bibliography{ajr}

\end{document}